\documentclass[conference,a4paper]{IEEEtran}
\IEEEoverridecommandlockouts
\usepackage{cite}
\usepackage{amsmath,amssymb,amsfonts}
\usepackage{algorithmic}
\usepackage{algorithm}
\usepackage{graphicx}
\usepackage{textcomp}
\usepackage{multirow}
\usepackage{booktabs}
\usepackage{color}
\usepackage[T1]{fontenc}
\usepackage[table,xcdraw]{xcolor}
\def\BibTeX{{\rm B\kern-.05em{\sc i\kern-.025em b}\kern-.08em
T\kern-.1667em\lower.7ex\hbox{E}\kern-.125emX}}
\begin{document}

\title{Approximating Energy-Regulation Feasible Region of Virtual Power Plants: A Data-driven Inverse Optimization Approach
}

\author{\IEEEauthorblockN{Ruike Lyu$^{1}$, Hongye Guo$^{1}$, and Qixin Chen*$^{1}$}
  \IEEEauthorblockA{Department of Electrical Engineering, Tsinghua University, Beijing, China$^{1}$}
  \IEEEauthorblockA{Email: qxchen@tsinghua.edu.cn$^{1}$}

  \thanks{This work was supported by National Natural Science Foundation of China under Grant 52107102 and the Major Smart Grid Joint Project of the National Natural Science Foundation of China and State Grid under Grant U2066205. Accepted by PESGM2024.}

}

\maketitle

\begin{abstract}
  System operators will probably allow virtual power plants (VPPs) to submit their feasible region (FR) for market clearing and dispatch. A VPP needs to determine its FR to submit as a whole based on the individual operation model of its internal distributed energy resources (DERs), which is an FR aggregation problem. Existing FR aggregation approaches rely on analytical methods, which have issues with generality and adaptability. In this paper, we propose a data-driven approach to approximate the energy-regulation FR of VPPs. It adopts the virtual battery model to approximate the aggregate FR of a VPP and determines the model parameters through inverse optimization based on generated multi-scenario operation data using the original operation model. Numerical tests verified the accuracy of the proposed method. We believe that our work helps to better leverage the flexibility of DERs.
\end{abstract}


\section{Introduction}

To better integrate the emerging distributed energy resources (DERs) such as energy storage in power systems, independent system operators (ISOs) are introducing new bidding models to represent their feasible regions (FRs), which are different from conventional units~\cite{order_841}. For example, new bidding models allow energy storage to submit their charging/discharging power limits and energy constraints so that their charging/discharging schedules can be determined through dispatch/market clearing. The rationale behind this change is that by incorporating the FR of DERs into the dispatch model, the resource mismatch issue that may arise from self-dispatch (DERs deciding their scheduling by themselves) can be avoided.

In the future, other types of resources may be allowed to submit their FRs in a similar way. Nevertheless, due to the typically small capacity of individual DERs, they usually need to be aggregated through a platform called a Virtual Power Plant (VPP) to reach the market threshold for participation~\cite{wei_bi-level_2018}. This means that a VPP needs to submit its FR as a whole through the bidding model specified by the ISO, which is more challenging than bidding for individual resources. The reason is that the process of obtaining the FR as a whole (i.e., the aggregate FR) based on the FRs of multiple resources is not a simple algebraic summation of the individual model parameters, but rather the addition of multiple sets (known as the Minkowski sum), for which there is currently no effective general solution method~\cite{wen_aggregate_2022}. Therefore, the common approach is to obtain a computationally feasible result using methods like inner or outer approximation~\cite{muller_aggregation_2019}, such as approximating the FR of DERs with ellipses to accelerate subsequent aggregation~\cite{chen_leveraging_2021}.

Currently, most aggregate FR approximation approaches are based on analytical methods, with limitations in terms of generality and adaptability in the context of VPP bidding. Firstly, FR aggregation methods usually rely on certain assumptions on the operation of individual DERs, for example, assuming that one DER's FR is independent of others~\cite{xu_hierarchical_2016}. However, in the scenario of VPP aggregating different types of resources participating in energy and regulation service markets, the energy-regulation FRs of the resources may be coupled~\cite{he_cooperation_2017}, so most FR aggregation methods cannot be directly applied. Secondly, typical aggregate FR approximation methods are tied with the form of individual FRs, such as box~\cite{chen_aggregate_2020} or ellipse~\cite{chen_leveraging_2021} as mentioned earlier, but these FR forms may not be consistent with the bidding model specified by ISO (e.g., battery models). As a result, analytical methods for aggregate FR approximation cannot be adaptively applied in VPP bidding.

Some recent studies have attempted to adopt FRs in the form of virtual battery models to address the above-mentioned issues, that is, approximating FR using the operating constraints of the battery, such as in~\cite{wang_aggregate_2021} and ~\cite{tan_optimal_2023}. Battery models have a more intuitive physical meaning and are consistent with the current bidding models. However, in these studies, the determination of FR parameters still relies on the analytical model of the DERs, making it hard to adapt to new types of resources. Another approach is to apply data-driven methods: the authors of Ref.~\cite{taheri_data-Driven_2022} propose using the original DER model to label randomly generated operation data as feasible or infeasible, which are used to train a convex quadratic classifier to approximate the aggregate FR. However, considering the temporal coupling characteristics of the DERs' operation, the required training data will exponentially increase with the number of time periods, and the trained classifier needs to be converted into the required form of the bidding model, introducing extra errors.

In summary, although data-driven methods may address the generality and adaptability issues of analytical FR aggregation methods, the relevant research is still in its early stages. Furthermore, most FR-aggregation studies focus on the FR for energy dispatch, while for resources such as energy storage that can quickly adjust their output, jointly participating in energy and regulation markets is more profitable~\cite{he_optimal_2016}. However, few studies have discussed the aggregation of energy-regulation FR of VPPs~\cite{yi_aggregate_2021}.

In this paper, we propose a data-driven approach for approximating the aggregate energy-regulation FR of VPPs, which adopts the virtual battery model as the form of the aggregate FR and fits the model parameters based on inverse optimization using generated operation data.
The basic idea is to first generate optimal dispatch results in different scenarios using the original DER operating models. These dispatch results are regarded as the optimal solutions based on the aggregate FR model in the corresponding scenarios, which are used as inputs for the inverse optimization problem to fit the parameters of the aggregate FR model.
Compared to existing methods, our approach combines the advantages of virtual battery models and data-driven methods: it can directly determine the FR parameters for bidding without relying on specific DER operating models.

\section{Problem Description}\label{sec_problem_description}

\begin{figure}[!t]
  \centering
  \includegraphics[width=1.5in]{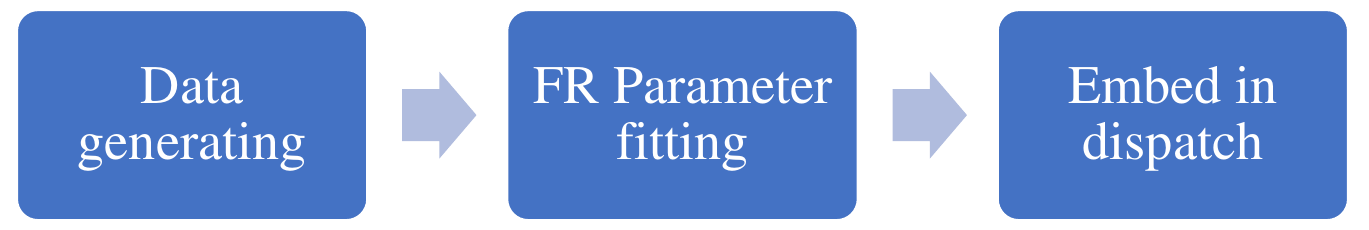}
  \caption{Framework of the proposed method: First, the original operation model of the VPP is used to generate the optimal dispatch results under different scenarios. Then, the dispatch results are used to fit the parameters of the approximate FR model based on inverse optimization. The obtained model can be used to approximate the aggregated FR for dispatching.}
  \label{fig_framework}
\end{figure}

Consider a VPP that aggregates DERs to jointly provide energy and frequency regulation services to the power system. According to the bidding model specified by the ISO, the VPP submits its FR for providing energy and regulation. Generally speaking, bidding models consist of two parts: FR and price. Nevertheless, here we focus on the FR part, and the treatment of price will be briefly discussed in the following sections
Based on the submitted FR and the bids/offers of other units and loads, the ISO solves the dispatching problem (or the clearing problem in a market environment) to determine the VPP's hourly energy output and regulation capacity. The VPP then schedules its internal DERs to meet the dispatch results as much as possible, as deviating from the dispatch results may incur additional penalties.

Without loss of generality, the FR $\mathcal{P}$ of the VPP is energy-regulation coupled, represented formally as:
\begin{align}\label{abstract_fr}
  \mathcal{P} := \{ & \ (p_{[T]}, r_{[T]}) : \exists y_{[I],[T]},                             \\
  \nonumber         & \ {\rm s.t.} \  h(p_{[T]}, r_{[T]}, y_{[I],[T]}; \theta_{[I]}) \le 0\}.
\end{align}
where $T$ and $I$ represent the number of scheduling time intervals and DERs, respectively. The variables with subscript $[T]$ is a shorthand notation for a T-dimensional vector, for example, $p_{[T]} = (p_1, p_2, ..., p_T)$ and $r_{[T]} = (r_1, r_2, ..., r_T)$ represent the hourly (baseline energy) output and regulation capacity of the VPP, respectively.
$h(\cdot)$ represents the original operational constraints of the VPP for providing energy-frequency services, with $y_{[I][T]}$ and $\theta_{[I]}$ representing the variables and parameters for the operation of the DERs within the VPP, respectively.

In Equation (\ref{abstract_fr}), $h(\cdot)$ includes the power and energy constraints of thousands of DERs within the VPP. The form of $\mathcal{P}$ for submission consists of a set of linear constraints on $(p, r)$, formally represented as $A(p, r) \le b$.
The calculation of $\mathcal{P}$ is a projection problem to eliminate the internal variables $y$ from the original constraints $h(\cdot)$, which is an NP-hard problem. Computing an exact $\mathcal{P}$ is often infeasible when $I$ is in the thousands. Instead, our main idea is to first approximate the original constraints $h(\cdot)$ with $h^{\rm a}(\cdot)$, the constraints of an approximate model with a much smaller dimension, and then construct an approximate FR $\mathcal{P}^{\rm a}$ based on $h^{\rm a}(\cdot)$:
\begin{align}\label{abstract_fr_surrogate}
  \mathcal{P}^{\rm a} := \{ & \ (p_{[T]}, r_{[T]}) : \exists y^{\rm a}_{[I^{\rm a}],[T]},                                                     \\
  \nonumber                 & \ {\rm s.t.} \  h^{\rm a}(p_{[T]}, r_{[T]}, y^{\rm a}_{[I^{\rm a}],[T]}; \theta^{\rm a}_{[I^{\rm a}]}) \le 0\}.
\end{align}
where $I^{\rm a}$ is the number of DERs in the approximate model. If $I^{\rm a} << I$, then using the constraints of the approximate model to derive $\mathcal{P}^{\rm a}$ or directly submitting $h^{\rm a}(\cdot)$ will be much less computationally intensive. For example, for a VPP consisting of $I = 1000$ electric vehicles, we can approximate the original constraints using an approximate model with $I^{\rm a} = 2$. However, it is necessary to carefully select the parameters $\theta^{\rm a}$ of the approximate model (which is the parameters of the battery model in this paper) to ensure that $\mathcal{P}^{\rm a}$ is sufficiently close to $\mathcal{P}$ (Fig.~\ref{fig_framework}), which is the main issue we focus on, represented as:
\begin{equation}\label{abstract_fr_approximation}
  \underset{\theta^{\rm a}}{\rm min.} J(\mathcal{P}, \mathcal{P}^{\rm a}(\theta^{\rm a}))
\end{equation}
where loss function $J(\cdot, \cdot)$ measures the distance between $\mathcal{P}$ and $\mathcal{P}^{\rm a}$.

\section{Methodology}\label{sec_method}

\subsection{Formulation of the approximate model}

The virtual battery (VB) model describes the operation of DERs using time-coupled power and energy constraints, which can be applied to model common resources such as energy storages, electric vehicles and thermostatically controlled loads. We employ the VB model to formulate the coupled energy-regulation constraints of the DERs. For resource $i$, its hourly output needs to satisfy the operation constraints for power limits (\ref{model_energy_powerLimit}), energy limits (\ref{model_energy_energyLimit}), change of energy (\ref{model_std_energyChange}), and initial energy (\ref{model_energy_energyInit}):
\begin{subequations}\label{model_energy}
  \begin{equation}\label{model_energy_powerLimit}
    0 \le p^{\rm dis(ch)}_{t, i} \le \overline{p}^{\rm dis(ch)}_{t, i} \ : \ \underline{\mu}^{\rm pd(c)}_{t, i}, \overline{\mu}^{\rm pd(c)}_{t, i}, \ \forall  t
  \end{equation}\vspace{-4ex}
  \begin{equation}\label{model_energy_energyLimit}
    \underline{e}_{t, i} \le e_{t, i} \le \overline{e}_{t, i} \ : \ \underline{\mu}^{\rm e}_{t, i}, \overline{\mu}^{\rm e}_{t, i},\ \forall t
  \end{equation}\vspace{-4ex}
  \begin{align}\label{model_std_energyChange}
    e_{t+1, i} = & e_{t, i} +
    (\eta^{\rm ch} p^{\rm ch}_{t, i} -  \frac{1}{\eta^{\rm dis}} p^{\rm dis}_{t, i}) \Delta t \ : \ \lambda^{\rm e}_{t, i},\ \forall t
  \end{align}\vspace{-4ex}
  \begin{equation}\label{model_energy_energyInit}
    e_{t, i} = e_{0, i} \ : \ \lambda^{\rm e0}_{i},\ t = 1
  \end{equation}
\end{subequations}
where $p_{t, i}$ and $r_{t, i}$ represent the output and regulation capacity for time interval $t$; $e_{t, i}$ is the energy at the end of $t$; the underlined/overlined parameters represents the lower/upper limit of the corresponding variables; the superscript dis(d)/ch(c) represent discharging to/charging from the grid; $e_{0, i}$ is the initial energy. $\eta$ is the (dis)charging efficiency. $\mu$ and $\lambda$ are the Lagrangian multipliers of the corresponding constraints.

The constraints for resource $i$ providing regulation include non-negative capacity (\ref{model_reg_nonnegative}), power capacity limits (\ref{model_reg_powerLimit1}-\ref{model_reg_powerLimit2}), and maintenance time requirements (\ref{model_reg_req1}-\ref{model_reg_req2}):
\begin{subequations}\label{model_reg}
  \begin{equation}\label{model_reg_nonnegative}
    r_{t, i} \ge 0 \ : \ \underline{\mu}^{r}_{t, i}, \ \forall t
  \end{equation}\vspace{-3ex}
  \begin{equation}\label{model_reg_powerLimit1}
    p^{\rm dis}_{t, i} - p^{\rm ch}_{t, i} + r_{t, i} - \overline{p}^{\rm dis}_{t, i} \le 0
    \ : \ \mu^{rpd}_{t, i}, \ \forall t
  \end{equation}\vspace{-3ex}
  \begin{equation}\label{model_reg_powerLimit2}
    - p^{\rm dis}_{t, i} + p^{\rm ch}_{t, i} + r_{t, i} - \overline{p}^{\rm ch}_{t, i} \le 0
    \ : \ \mu^{rpc}_{t, i}, \ \forall t
  \end{equation}\vspace{-3ex}
  \begin{equation}\label{model_reg_req1}
    \eta^{\rm ch} (r_t - p^{\rm dis}_{t, i} + p^{\rm ch}_{t, i}) \Delta t^{req} + e_{t, i} - \overline{e}_{t, i} \le 0
    \ : \ \mu^{rec}_{t, i}, \ \forall t
  \end{equation}\vspace{-3ex}
  \begin{equation}\label{model_reg_req2}
    \frac{1}{\eta^{\rm dis}} (r_t + p^{\rm dis}_{t, i} - p^{\rm ch}_{t, i}) \Delta t^{req} - e_{t, i} + \underline{e}_{t, i} \le 0
    \ : \ \mu^{red}_{t, i}, \ \forall t
  \end{equation}
\end{subequations}
where $\Delta t^{req}$ represents the duration required for the resource to maintain its maximum regulation output, for example, 15 minutes in PJM~\cite{he_optimal_2016}. The bids of the VPP are the aggregated bids of the DERs, which are represented as:
\begin{equation}\label{model_vpp_aggBid}
  p_t = \sum_{i=1}^{I} (p^{\rm dis}_{t, i} - p^{\rm ch}_{t, i}), \ r_t = \sum_{i=1}^{I} r_{t, i}, \ \forall t
\end{equation}\vspace{-3ex}

In the above model, the parameters are $\theta = \{\overline{p}^{\rm dis(ch)}_{t, i}, \underline{e}_{t, i}, \overline{e}_{t, i}, e_{0, i}, \forall t \forall i\}$, and the internal variables are $y = \{p_{t, i}, r_{t, i}, e_{t, i}, \forall t \forall i\}$. The constraints $h(\cdot)$ are the combination of Equations (\ref{model_energy})-(\ref{model_vpp_aggBid}).

Note that although we omitted the superscript a representing the approximate model above, which implies the assumption of formal consistency between the approximate model and the original model, our method does not depend on this assumption. In the case study, we tested the generality of the proposed method in scenarios where the approximate model and the original model are different in formulations.

\subsection{Determining the parameters of the approximate model}

The parameters of the approximate model $\theta^{\rm a}$ need to be determined based on the original operational model of the VPP. Here, we propose using data-driven inverse optimization to determine $\theta^{\rm a}$, which consists of two steps: data generation and parameter fitting. The basic logic is that if the dispatching results given by the ISO based on the $\mathcal{P}^{\rm a}(\theta^{\rm a})$ submitted by the VPP can approximate the dispatching results based on $\mathcal{P}$ under various dispatching scenarios, then $\mathcal{P}^{\rm a}$ can be considered sufficiently close to $\mathcal{P}$.

\subsubsection{Data generation}
In theory, to obtain the dispatching results after submitting the original FR, it is necessary to model all units and loads within the scope of ISO. However, in practice, the capacity of a VPP is very small, so it can be modeled as a price taker. Thus, under the assumption of perfect competition, the dispatching scenarios can be represented by energy-regulation price scenarios~\cite{zhou_incentive-compatible_2023}. Thus, the joint energy-regulation dispatching of the VPP is equivalent to the optimization problem below:
\begin{subequations}\label{model_bid_primal}
  \begin{equation}\label{model_vpp_profit}
    \underset{p_{[T]}, r_{[T]}}{\rm min.} - \sum_{t=1}^{T}(Pr^{\rm e}_{t} p_{t} + Pr^{\rm r}_{t} r_{t}
    - Pr^{\rm deg}(p_{t}))\Delta t
  \end{equation}\vspace{-2ex}
  \begin{equation}\label{model_vpp_constraint}
    {\rm s.t.} \ (p_{[T]}, r_{[T]}) \in \mathcal{P}
  \end{equation}
\end{subequations}
where $Pr^{\rm e}_{t}$ is the energy price; $Pr^{\rm r}_{t} = s^{\rm perf} (Pr^{\rm cap}_{t} + Pr^{\rm mil}_{t} a^{\rm mil}_{t})$ is the equivalent regulation capacity price for the VPP; $s^{\rm perf}$ is the performance score, $Pr^{\rm cap(mil)}_{t}$ is the regulation capacity(mileage) price, and $a^{\rm mil}_{t}$ is the expected regulation mileage for $t$. $Pr^{\rm deg}$ is the cost function of the resources, which can be included in the bidding model along with $\mathcal{P}$, but its specific determination is not the focus of this paper. In the following text, we assume that $Pr^{\rm deg}$ is directly proportional to the total discharge power of the VPP.

Given the energy-regulation price $(Pr^{\rm e}_{[T]}, Pr^{\rm r}_{[T]})$, solving (\ref{model_bid_primal}) using the original VPP model can obtain the corresponding scheduling results $(p_{[T]}, r_{[T]})$, thus generating a price-dispatch dataset $D = \{(Pr^{\rm e(k)}_{[T]}, Pr^{\rm r(k)}_{[T]}, p^{(k)}_{[T]}, r^{(k)}_{[T]})\}$, where $k$ and $K$ represent the index and the total number of scenarios, respectively.

In practice, the original constraints $h(\cdot)$ can be directly substituted to replace (\ref{model_vpp_constraint}), as the optimal solution of the problem remains unchanged before and after projection~\cite{tan_non-iterative_2023}.

\subsubsection{Parameter fitting}
The idea is to treat the generated scheduling results in data set $D$ as the optimal solutions (with some noise added) of the optimization problem (\ref{model_bid_primal}) formulated by the approximate model, and to fit the problem parameters using $D$. This can be represented by the following data-driven inverse optimization (DDIO) problem:
\begin{subequations}\label{model_inverse_fit}
  \begin{equation}\label{model_inverse_fit_lossFunction}
    \underset{\theta^{\rm a}}{\rm min.} J = \frac{1}{K} \sum_{k = 1}^{K} (
    ||p^{\rm a(k)}_{[T]} - p^{(k)}_{[T]}||^2 + ||r^{\rm a(k)}_{[T]} - r^{(k)}_{[T]}||^2
    )
  \end{equation}\vspace{-1ex}
  \begin{equation}\label{model_inverse_fit_constraint_1}
    {\rm s.t.} \ {\rm strong \ duality:} (\ref{model_vpp_profit}) = g(\mu^{(k)}, \lambda^{(k)}), \forall k.
  \end{equation}\vspace{-3ex}
  \begin{equation}\label{model_inverse_fit_constraint_2}
    {\rm primal \ feasibility:} h^{\rm a}(p^{\rm a(k)}_{[T]}, r^{\rm a(k)}_{[T]}, y^{\rm a(k)}; \theta^{\rm a}) \le 0, \forall k.
  \end{equation}\vspace{-3ex}
  \begin{equation}\label{model_inverse_fit_constraint_3}
    {\rm dual \ feasibility:} (\ref{model_dual_feasibility}); {\rm stationarity:} (\ref{model_coefficients}), \forall k.
  \end{equation}
\end{subequations}
where $g(\mu^{(k)}, \lambda^{(k)})$ is the Lagrangian dual function of (\ref{model_bid_primal}); (\ref{model_inverse_fit_constraint_1}-\ref{model_inverse_fit_constraint_3}) are the strong duality, primal feasibility, dual feasibility, and stationarity conditions, respectively, guaranteeing the optimality of $(p^{\rm a(k)}_{[T]}, r^{\rm a(k)}_{[T]})$. The Lagrangian dual function and the corresponding constraints are given in Appendix~\ref{app_dual}.

\subsection{An iterative algorithm for solving the DDIO problem}

Determining the optimal parameters (\ref{model_inverse_fit}) is a high-dimensional nonlinear optimization problem, especially challenging to solve when $K$ is large. One approach is to update $\theta^{\rm a}$ step by step based on the structure of the problem using Newton's method~\cite{tan_data-driven_2023}, but this relies on specific forms of the optimization problem. Instead, we adopt a commonly used method in machine learning, where we update the parameters iteratively using the data of one scenario at a time. When updating the parameters with data of scenario $k$, solve the following problem:
\begin{subequations}\label{model_inverse_fit_update}
  \begin{equation}\label{model_inverse_fit_update_lossFunction}
    \underset{\theta^{\rm a(k)}}{\rm min.} (
    ||p^{\rm a(k)}_{[T]} - p^{(k)}_{[T]}||^2 + ||r^{\rm a(k)}_{[T]} - r^{(k)}_{[T]}||^2
    ) + \alpha ||\theta^{\rm a(k)} - \theta^{\rm a}||^2
  \end{equation}\vspace{-3ex}
  \begin{equation}\label{model_inverse_fit_update_constraint}
    {\rm s.t.} (\ref{model_inverse_fit_constraint_1}-\ref{model_inverse_fit_constraint_3}), {\rm for} \ k.
  \end{equation}
\end{subequations}
where $\theta^{\rm a}$ is the parameter value of the previous iteration. The last term in (\ref{model_inverse_fit_update_lossFunction}) is added to help converge, where $\alpha$ adaptively increases as with the iterations. (\ref{model_inverse_fit_update}) can be directly fed to commercial solvers. In practice, a maximum computation time can be set for each iteration, and the obtained feasible solution when at the maximum time can be used for update. This strategy can easily be implemented by commercial solvers. The overall iterative procedure is given in Algorithm~\ref{alg_framework}. For initializing $\theta^{\rm a}$, we summed the original model parameters directly and distributed them equally to the approximate model in numerical tests.

\begin{algorithm}[!t]
  \caption{Iterative solution of the DDIO problem.}
  \label{alg_framework}
  \begin{algorithmic}[1]
    \renewcommand{\algorithmicrequire}{\textbf{Input:}}
    \REQUIRE
    maximum number of iterations $N$, price-dispatch data set $\{(Pr^{\rm e(k)}_{[T]}, Pr^{\rm r(k)}_{[T]}, p^{(k)}_{[T]}, r^{(k)}_{[T]})\}$.
    \renewcommand{\algorithmicrequire}{\textbf{Output:}}
    \REQUIRE an approximate model parameterized by $\theta^{\rm a}$.
    \STATE initialize the parameters $\theta^{\rm a}$, set $n = 0$.
    \WHILE{$n \le N$}
    \STATE let $k = n \ {\rm mod} \ K + 1$, solve (\ref{model_inverse_fit_update}) to obtain $\theta^{\rm a(n)}$.
    \STATE update $\theta^{\rm a}$ by the average of the last min$\{K, n\}$ $\theta^{\rm a(n)}$s
    \STATE \textbf{if} $\theta^{\rm a}$ converges, \textbf{break}.
    \STATE set $n = n + 1$.
    \ENDWHILE
  \end{algorithmic}
\end{algorithm}

\section{Numerical Results}\label{sec_numerical}

We used Gurobi (V10.0.0) with YALMIP~\cite{Lofberg2004} to solve the optimization problems on a workstation with an Intel Core i9-10900X CPU (3.7 GHz) and 128 GB RAM.

\subsection{Scenario setting and data generation}

\begin{figure}[!t]
  \centering
  \includegraphics[width=1.5in]{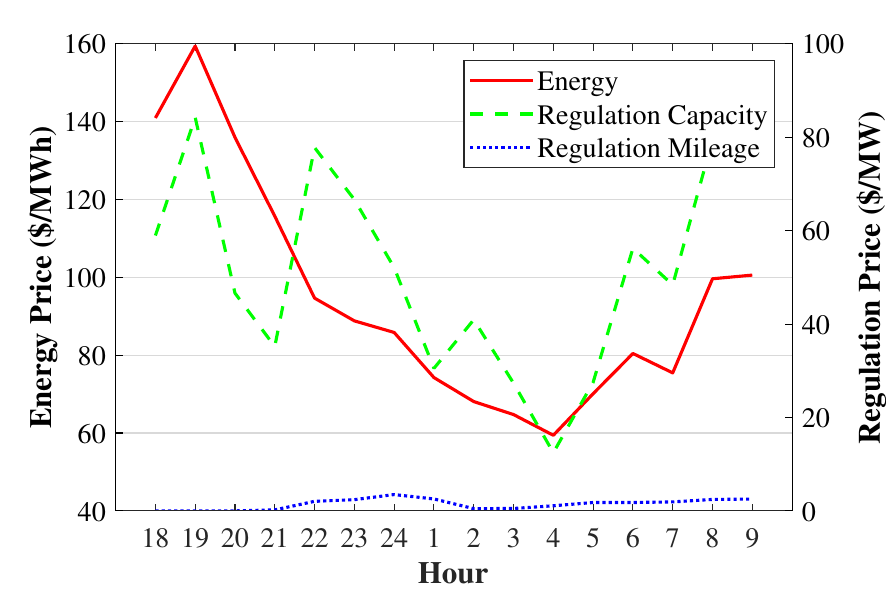}
  \includegraphics[width=1.5in]{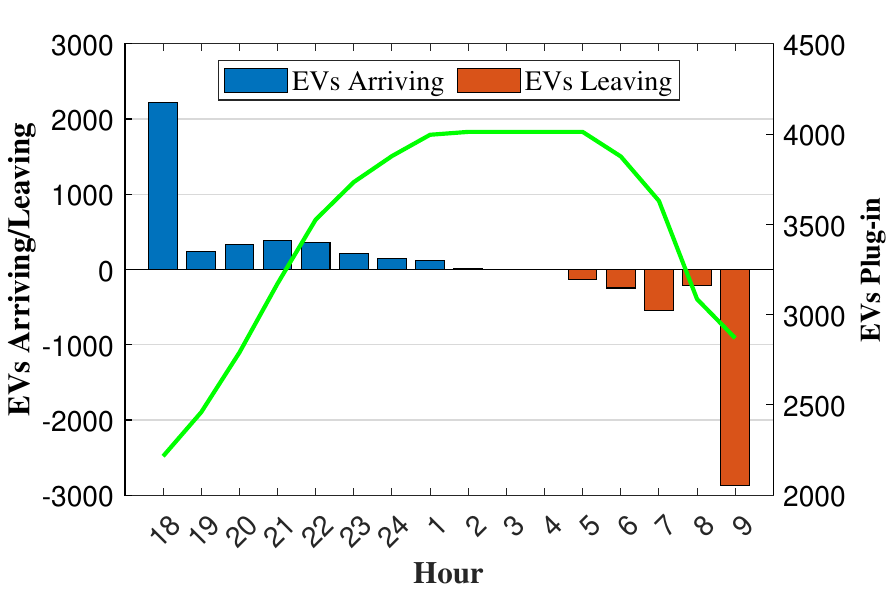}
  \caption{Scenario setting: (left) typical energy / regulation prices of PJM (July 27, 2022). (right) The arriving / leaving time of the EVs.}
  \label{fig_scenario}
\end{figure}

We borrowed the scene of aggregating 4000 electric vehicles (EVs) for providing energy and regulation services from Ref.~\cite{lyu_co-optimizing_2023}. The EV batteries were charged starting from 0kWh using bidirectional 7.68kW charging piles, and it was required that the battery level at departure was not less than 80\% of the maximum capacity of 50kWh. The compensation for battery discharging was set at \$0.1/kWh. Readers can refer to Ref.~\cite{lyu_co-optimizing_2023} for detailed scenario setting.

We used historical energy and regulation market prices from PJM in July 2022 to represent different dispatching scenarios. For the scenarios from July 1st to 20th, the original model of the 4000 EVs was used to generate the energy-regulation dispatching results ($2 \times T$ matrices for each scenario, $T = 16$), which were used for fitting the aggregated FR parameters. The data from July 21st to 30th was used to test the error between the dispatching results using the fitted FR and the results using the original model. Other market parameters include the length of (scheduling) time interval $\Delta t = 1$ hour,
the required maintenance time $\Delta t^{\rm req} = 0.25$ hours, and
the performance score $s^{\rm perf} = 0.984$. The regulation mileage was estimated using historical data.

To test the generality of the proposed method, which means it can be applied to aggregate different forms of individual FRs, we used two DER operational models for providing regulation. The first model is the commonly used \textbf{decoupled operation} model, assuming each DER deploys regulation according to its own capacity, so the FR of each DER is independent of others~\cite{vagropoulos_optimal_2013}. The second model optimizes the power allocation among the DERs, resulting in \textbf{coupled operation} model~\cite{lyu_co-optimizing_2023}. To compare the effectiveness of the proposed method, we also tested the results using the \textbf{outer approximation} method from Ref.~\cite{xu_hierarchical_2016}.

\subsection{Results and comparison}

\begin{figure}[!t]
  \centering
  \includegraphics[width=2.0in]{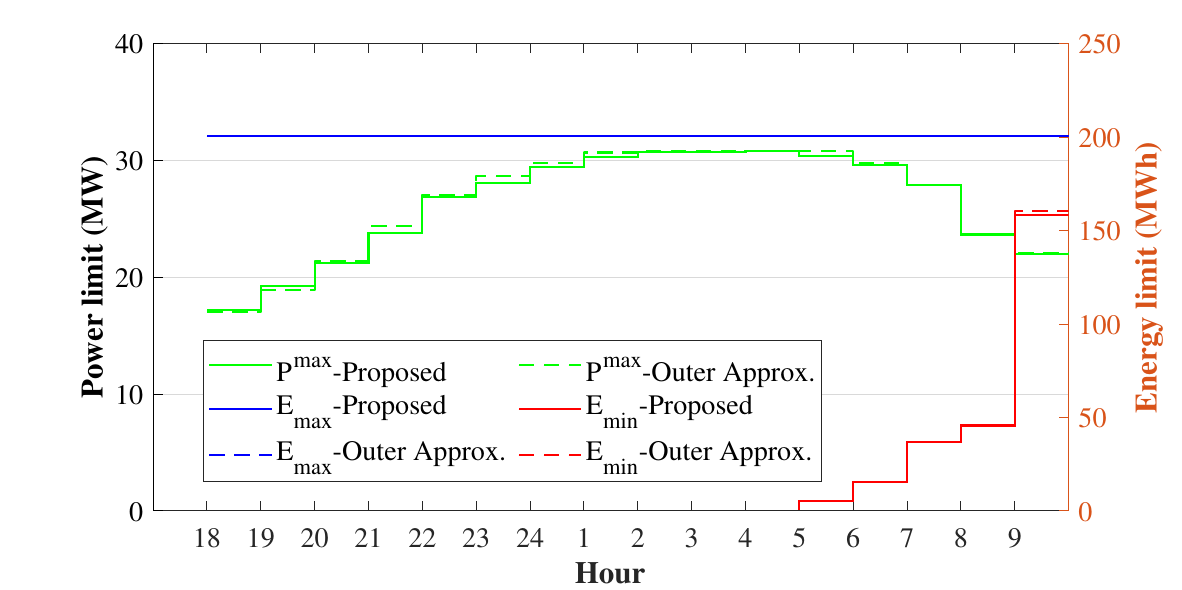}
  \caption{Aggregated FR parameters of the decoupled operation model given by the proposed method (I = 1) / outer approximation.}
  \label{fig_parameter}
\end{figure}

\begin{figure}[!t]
  \centering
  \includegraphics[width=1.5in]{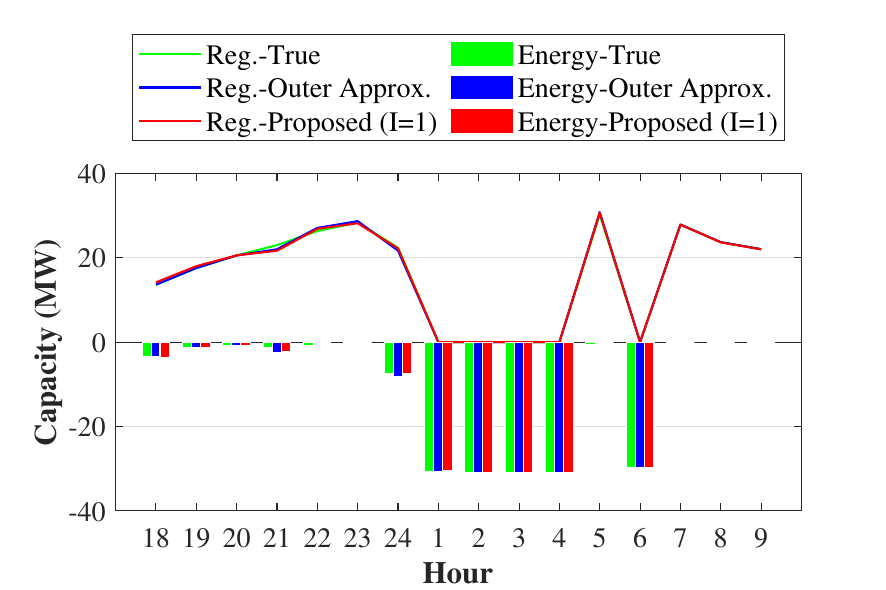}
  \includegraphics[width=1.5in]{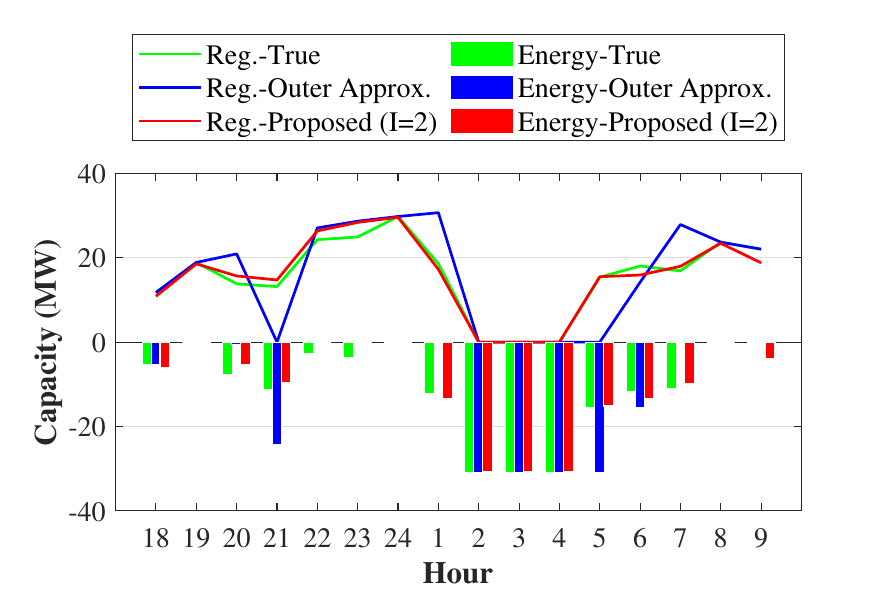}
  \caption{Typical dispatch results based on the aggregated FR (July 27, 2022): (left) Decoupled DER operation model; (right) Coupled DER operation model.}
  \label{fig_dispatch}
\end{figure}

\begin{table}[!t]\renewcommand{\arraystretch}{1}
  \caption{Normalized Error (MAE) of Aggregate FR}
  \label{tab_mae}
  \centering
  \begin{tabular}{cccc}
    \toprule
    \begin{tabular}[c]{@{}c@{}}DER operation \\ model\end{tabular} & \begin{tabular}[c]{@{}c@{}}Outer \\ Approx.~\cite{xu_hierarchical_2016}\end{tabular} & \begin{tabular}[c]{@{}c@{}}Proposed \\ (I=1)\end{tabular} & \begin{tabular}[c]{@{}c@{}}Proposed \\ (I=2)\end{tabular} \\ \hline
    \begin{tabular}[c]{@{}c@{}}Decoupled \\ operation~\cite{vagropoulos_optimal_2013}\end{tabular} & 1.8\%                      & \textbf{1.6\%}             & 1.7\%                      \\ \hline
    \begin{tabular}[c]{@{}c@{}}Coupled \\ operation~\cite{lyu_co-optimizing_2023}\end{tabular} & 13.6\%                     & 13.4\%                     & \textbf{7.5\%}             \\ \hline
  \end{tabular}
\end{table}

The typical aggregate FR parameters obtained through the proposed method and the outer approximation method are given in Fig.~\ref{fig_parameter}. The results from both methods are close. The obtained aggregate FRs were then embedded in the dispatch problem, and the dispatching results of the VPP under different scenarios are compared with the results obtained using the original model, with the typical results shown in Fig.~\ref{fig_dispatch} and the errors presented in Table~\ref{tab_mae}. The errors were measured using the ratio of mean absolute error (MAE) to the maximum true value (results using the original model). The results show that when the decoupled operation model is used as the original model, the errors of the compared aggregate FR approximation methods are small (less than 2\%). However, with coupling in the DERs' FRs, the errors of the outer approximation method exceeded 10\%, while the proposed method reduced the errors by over 40\% compared to the outer approximation method. This indicates that the proposed method is at least not inferior to typical FR approximation methods in terms of accuracy in the tests and, due to its better generality, performs better in some scenarios.

\section{Conclusion}\label{sec_conclusion}

The VPP could be allowed to submit its aggregate FR to the ISO so that the ISO can conduct dispatch within its operating boundaries. To address the generality and adaptability issues in the analytical method of FR aggregation, this paper proposes a data-driven method that fits the VPP's aggregate FR parameters based on multi-scenario optimal dispatch data generated by the original operating model, thereby avoiding bottom-up analytical analysis and aggregation. Compared with existing approximation methods for FR aggregation, our method can directly determine the FR parameters to submit to the ISO without relying on the analytical form of the DER operating model. Numerical results validated the advantages of the proposed method compared to existing methods. In cases where there are coupled operating constraints for individual DERs, the aggregate FR error was reduced by over 40\%.

\appendices
\section{Dual function and constraints of problem (\ref{model_bid_primal})}\label{app_dual}
We omit the superscript $k$ here. $\forall k$, we have: Lagrangian dual function:
\begin{equation}\label{model_dual_function}
  \begin{aligned}
    g(\mu, \lambda) = & \sum_{i=1}^{I} (\sum_{t=1}^{T}(
    \underline{e}_{t, i} \underline{\mu}^{\rm e}_{t, i} - \overline{e}_{t, i} \overline{\mu}^{\rm e}_{t, i}
    - \overline{p}^{\rm dis}_{t, i} \mu^{\rm rpd}_{t, i}
    - \overline{p}^{\rm ch}_{t, i} \mu^{\rm rpc}_{t, i}            \\
                      & - \overline{e}_{t, i} \mu^{\rm rec}_{t, i}
    + \underline{e}_{t, i} \mu^{\rm red}_{t, i}) + e_{0, i} \lambda^{\rm e0}_{i})
  \end{aligned}
\end{equation}\vspace{-3ex}

Dual feasibility ($\forall i$):
\begin{equation}\label{model_dual_feasibility}
  (\underline{\mu}_{t, i}^{\rm pd(c)}, \underline{\mu}_{t, i}^{\rm e}, \overline{\mu}_{t, i}^{\rm e}, \underline{\mu}_{t, i}^{\rm r}, \mu_{t, i}^{\rm r p d(c)}, \mu_{t, i}^{\rm r e d(c)}) \ge 0, \ \forall t, \forall i
\end{equation}\vspace{-3ex}

Stationarity ($\forall i$):
\begin{subequations}\label{model_coefficients}
  \begin{equation}\label{model_coefficients_pdis}
    \begin{aligned}
      Pr^{\rm deg} - Pr^{\rm e}_t - \underline{\mu}_{t, i}^{\rm p d}-\frac{\Delta t}{\eta^{\rm dis}} \lambda^{\rm e}_{t, i} + \mu_{t, i}^{\rm r p d} - \mu_{t, i}^{\rm r p c} \\
      + \frac{\Delta t^{\rm req}}{\eta^{\rm dis}} \mu_{t, i}^{\rm r e d} - \eta^{\rm ch} {\Delta t}^{\rm req} \mu_{t, i}^{\rm r e c}=0, \ \forall t
    \end{aligned}
  \end{equation}\vspace{-3ex}
  \begin{equation}\label{model_coefficients_pch}
    \begin{aligned}
      Pr^{\rm e}_{t} - \underline{\mu}_{t, i}^{\rm p c}+\Delta t \eta^{\rm ch} \lambda^{\rm e}_{t, i} - \mu_{t, i}^{\rm r p d} + \mu_{t, i}^{\rm r p c} \\
      - \frac{\Delta t^{\rm req}}{\eta^{\rm dis}} \mu_{t, i}^{\rm r e d} + \eta^{\rm ch} {\Delta t}^{\rm req} \mu_{t, i}^{\rm r e c} = 0, \ \forall t
    \end{aligned}
  \end{equation}\vspace{-3ex}
  \begin{equation}\label{model_coefficients_e0}
    \lambda_{t, i}^{\rm e} -\lambda^{\rm e0} + \mu_{t, i}^{\rm r e c} - \mu_{t, i}^{\rm r e d}=0, \ t = 1
  \end{equation}\vspace{-3ex}
  \begin{equation}\label{model_coefficients_e}
    -\underline{\mu}_{t, i}^{\rm e} + \overline{\mu}_{t, i}^{\rm e}+\lambda_{t, i}^{\rm e} - \lambda_{t - 1}^{\rm e}+\mu_{t, i}^{\rm r e c}-\mu_{t, i}^{\rm r e d}=0, \ 1 < t < T + 1
  \end{equation}\vspace{-3ex}
  \begin{equation}\label{model_coefficients_eEnd}
    -\underline{\mu}_{t, i}^{\rm e} + \overline{\mu}_{t, i}^{\rm e} - \lambda_{t - 1}^{\rm e}=0, \ t = T + 1
  \end{equation}\vspace{-3ex}
  \begin{equation}\label{model_coefficients_r}
    - Pr_{t}^{\rm r} - \underline{\mu}_{t, i}^{\rm r}+\mu_{t, i}^{\rm r p d}+\mu_{t, i}^{\rm r p c} + \eta^{\rm ch} {\Delta t}^{\rm req} \mu_{t, i}^{\rm r e c} + \frac{\Delta t^{\rm req}}{\eta^{\rm dis}} \mu_{t, i}^{\rm r e d}=0, \ \forall t
  \end{equation}
\end{subequations}

\bibliographystyle{IEEEtran}
\bibliography{reference}

\end{document}